\documentclass[a4paper,12pt]{article}
\usepackage{amsfonts, amsmath, amssymb, amscd, latexsym, graphicx}
\newtheorem{thm}{Theorem}[section]
\newtheorem{prop}[thm]{Proposition}
\newtheorem{lem}[thm]{Lemma}

\newtheorem{defn}[thm]{Definition}

\newcommand{\demo}{ {\it   Proof. }}

\title{\Large On subgroup conjugacy separability\\ in the class of virtually free groups}

\author{O.\ Bogopolski \\ \small{Institute of Mathematics of} \\ \small{Siberian Branch of Russian Academy
of Sciences,} \\ {\small Novosibirsk, Russia}\\ {\small and D\"{u}sseldorf University, Germany} \\ \small{e-mail:
Oleg$\_$Bogopolski@yahoo.com}\\ \\
\framebox[30mm][c]{F.\ Grunewald} \footnote{During the process of writing this paper Fritz Grunewald died.
} \\
{\small D\"{u}sseldorf University, Germany}}
\large

\begin{document}

\maketitle

\begin{abstract}

A group $G$ is called subgroup conjugacy separable (abbreviated as SCS), if
any two finitely generated and non-conjugate subgroups of
$G$ remain non-conjugate in some finite quotient of $G$.
We prove that the free groups and the fundamental groups of finite trees of finite groups
with some normalizer condition are SCS.
We also introduce the subgroup into-conjugacy separability property and prove
that the above groups have this property too.

\end{abstract}

\bigskip

\section{Introduction}


The subgroup conjugacy separability (see Definition~\ref{def2}) is a residual property of groups, which logically continues the following series of well
known properties of groups: the residual finiteness, the conjugacy separa\-bility, and the subgroup separability (LERF).
These properties help to solve some algorithmic problems in groups
and they are also important in the theory of 3-manifolds.

Recall that a group $G$ is called {\it separable} (or residually finite, abbreviated as {\rm RF}), if for any two elements $x\neq y\in G$, there exists a homomorphism $\phi$ from $G$ to a finite group $\overline{G}$ such that $\phi(x)\neq\phi(y)$ in $\overline{G}$.

Similarly, $G$ is called
{\it conjugacy separable} (abbreviated as {\rm CS}), if for any two non-conjugate elements $x,y\in G$, there
exists a homomorphism $\phi$ from $G$ to a finite group $\overline{G}$ such that $\phi(x)$ is not conjugate to
$\phi(y)$ in $\overline{G}$.

If in these definitions we replace the words ``elements'' by the words ``finitely generated subgroups'', we
obtain the following two definitions (the first one is well known, and the second is new).



\begin{defn}\label{def1}{\rm A group $G$ is
called {\it subgroup separable} (or LERF, for locally extended residually finite),
if for any two finitely generated subgroups $H_1\neq H_2\leqslant G$, there
exists a homomorphism $\phi$ from $G$ to a finite group $\overline{G}$ such that $\phi(H_1)\neq\phi(H_2)$ in $\overline{G}$.}
\end{defn}


Note that this definition is equivalent to the usual one: $G$ is called subgroup separable if
for any finitely generated subgroup $H$ of $G$ and for any element $g\in G\setminus H$, there exists a
subgroup $K$ of finite index in $G$ such that $H\leqslant K$, but $g\notin K$.

The LERF property is useful in 3-manifold topology: if $\pi_1(M^3)$ is a LERF group and $S\rightarrow M^3$
is an immersion of an incompressible surface, then there is an embedding $S\hookrightarrow \widetilde{M}^3$
in a finite cover $\widetilde{M}^3$ of $M^3$. For more information, see the inspiring paper of P.~Scott~\cite{Scott} and an overview of D.~Wise in~\cite{W2}. It would be interesting to find applications of the following property in
topology.

\begin{defn}\label{def2} {\rm A group $G$ is called
{\it subgroup conjugacy separable} (abbreviated as SCS), if for any two finitely generated non-conjugate subgroups $H_1,H_2\leqslant G$, there exists a homomorphism $\phi$ from $G$ to a finite group $\overline{G}$ such that $\phi(H_1)$ is not conjugate to $\phi(H_2)$ in $\overline{G}$.}
\end{defn}

A.I.~Mal'cev was the first, who noticed, that finitely presented residually finite (resp. conjugacy separable) groups have solvable word problem (resp. conjugacy problem)~\cite{Mal}. Arguing in a similar way, one can show that finitely presented LERF groups
have solvable membership problem and that finitely presented SCS groups have solvable conjugacy problem
for finitely generated subgroups. The last means, that there is an algorithm, which given a finitely presented SCS group $G=\langle X\,|\,R\rangle$ and two finite sets of elements $U=\{u_1,\dots,u_n\}$ and $V=\{v_1,\dots,v_m\}$, decides whether the subgroups $\langle U\rangle$ and $\langle V\rangle$ are conjugate in $G$.



Clearly, any group with the property CS, LERF, or SCS is residually finite.
We conjecture, that SCS does not imply CS or LERF and conversely.


There is a lot of papers devoted to the properties RF, CS, and LERF.
We cite here some positive results
about the CS and LERF properties. It would be interesting to establish, which of the listed below groups
have the SCS property.

\medskip

{\bf The conjugacy separability} was established for

\vspace*{-3mm}
\begin{enumerate}

\item[--] virtually polycyclic groups (V.~Remeslen\-nikov~\cite{Rem1} and E.~Formanek~\cite{For})\vspace*{-3mm}

\item[--] finitely generated virtually free groups (see J.L.~Dyer~\cite{D} combined with
P.F.~Stebe~\cite{Stebe1})\vspace*{-3mm}

\item[--] groups which can be obtained from these groups by repeatedly forming free products with cyclic amalgamations
(L.~Ribes, D.~Segal und P.~Zalesskii~\cite{RSZ})\vspace*{-3mm}

\item[--] virtually surface groups\footnote{For free products of two
 free groups, amalgamated over a cyclic group, in particular for surface groups see
the paper of J.L.~Dyer~\cite{D1}. For Fuchsian groups see the paper of B.~Fine and G.~Rosenberger~\cite{FR}.}
and the fundamental groups of Seifert 3-manifolds (A.~Martino~\cite{Martino})\vspace*{-3mm}

\item[--] fundamental groups of finite, 1-acylindrical graphs of
free groups with finitely generated edge groups (O.~Cotton-Barratt und H.~Wilton~\cite{CW})\vspace*{-3mm}

\item[--] virtually limit groups (S.~Chagas and P.~Zalesskii~\cite{CZ0,CZ2})\vspace*{-3mm}

\item[--] finitely presented residually free groups (S.~Chagas and P.~Zalesskii~\cite{CZ1})\vspace*{-3mm}

\item[--] right angled Artin groups and their finite index subgroups (A.~Minasyan~\cite{Min}).\vspace*{-3mm}

\item[--] free products of CS groups (P.F. Stebe~\cite{Stebe1} and V.N. Remeslennikov~\cite{Rem2})\vspace*{-3mm}

\item[--] non-uniform arithmetic lattices of $SL_2(\mathbb{C})$ and consequently the Bianchi groups
(S.~Chagas and P.~Zalesskii~\cite{CZ0}; see also the paper of I. Agol, D.D. Long, and A.W.~Reid~\cite{ALR})

\end{enumerate}

{\bf The subgroup separability} was established for
\vspace*{-3mm}
\begin{enumerate}

\item[--] polycyclic groups  (A.I.~Mal'cev~\cite{Mal})\vspace*{-3mm}

\item[--] free groups (M.~Hall~\cite{Hall}\vspace*{-3mm}

\item[--] surface groups (P.~Scott~\cite{Scott})\vspace*{-3mm}

\item[--] limit groups (H.~Wilton~\cite{Wilt})\vspace*{-3mm}

\item[--] free products of LERF groups (R.G.~Burns~\cite{Burns} and N.S.~Romanovskii~\cite{Rom})\vspace*{-3mm}

\item[--] free products of two free groups
amalgamated along a cyclic group (A.M.~Brunner, R.G.~Burns and D.~Solitar~\cite{BBS}; see also a generalization of M.~Tretkoff~\cite{T})\vspace*{-3mm}

\item[--] free products of a LERF group $G$ and a free group $F$ amalgamated along a maximal cyclic subgroup in $F$ (R.~Gitik~\cite{G}) (Note, that the free product of two LERF groups amalgamated along a cyclic subgroup is not necessarily a LERF group (see~\cite{Rips} and~\cite{GR1})\vspace*{-3mm}

\end{enumerate}

\medskip

If $G$ splits as a finite graph of free groups with cyclic edge groups, then $G$ is LERF if and only if
$G$ does not contain a non-trivial element $a$, such that $a^n$ is conjugate to $a^m$
for some $n\neq \pm m$ (D.~Wise~\cite{W1}).

In~\cite{MR}, V.~Metaftsis and E.~Raptis are proved that a right-angled Artin group $G$ with associated graph $\Gamma$ is
subgroup separable if and only if $\Gamma$ does not contain a subgraph homeomorphic
to either a square or a path of length three.

P.~Scott in~\cite{Scott} showed, that LERF is inherited by subgroups and finite extensions, in particular
it is invariant under commensurability. In contrast, CS is not invariant under commensurability:
in~\cite{MM}, A.~Martino and A.~Minasyan
constructed a finitely presented CS-group, which has an index 2 subgroup without the CS property.
An example of a finitely generated (but not finitely presented)
non-CS-group $G$, containing a CS-subgroup of
index 2, was constructed by A. Gorjaga in~\cite{Gor}.



\medskip

{\bf The subgroup conjugacy separability property.} It seems that the SCS property is much harder to establish than
CS and LERF. One of the reasons is that this property has no an evident reformulation in terms of
the profinite topology on $G$.

Recall that the {\it profinite topology} on a group $G$ is the topology,
having the family of all cosets of subgroups of finite index in $G$ as a base of open sets.
Clearly, a finitely generated group $G$ is residually finite, respectively conjugacy separable or subgroup separable, if and only if one-element subsets of $G$, respectively conjugacy classes of one-element subsets,
or finitely generated subgroups are closed in the profinite topology. We conjecture, that the subgroup conjugacy separability for $G$ is not the same as the closeness, in the profinite topology, of the union of the conjugacy classes of any finitely generated subgroup of $G$.

We know only one paper on SCS (without restrictions on subgroups): in~\cite{GS}, F.~Grune\-wald and D.~Segal proved that all virtually polycyclic groups are subgroup conjugacy separable (see also Theorem 7 in Chapter 4 of~\cite{Segal}).

In this paper we consider finitely generated virtually free groups. These groups are subgroup separable (since they are commensurable with free groups) and they are
conjugacy separable (J.L.~Dyer~\cite{D} and  P.F.~Stebe~\cite{Stebe1}).
Therefore it is natural to ask, whether all finitely generated virtually free groups
are subgroup conjugacy separable.

Recall that every finitely generated virtually free group
is the fundamental group of a finite graph of finite groups (A.~Karrass, A.~Pietrowski and D.~Solitar~\cite{KPS}).
The main results of this paper are Theorems~\ref{main0},~\ref{main1}, and~\ref{main3},~\ref{main2}.

\begin{thm}~\label{main0} Free groups are subgroup
conjugacy separable.
\end{thm}

In the following definition we use the notations of Section~\ref{1}.

\begin{defn} {\rm We say that a finite graph of finite groups $(\mathbb{G},\Gamma)$ (and its fundamental group) satisfies the {\it normalizer condition}, if $|N_G(E):E|<\infty$ for each nontrivial subgroup $E$ of every edge group of $G=\pi_1(\mathbb{G},\Gamma)$.}
\end{defn}

For instance, $A\ast_C B$ satisfies the normalizer condition, if $A,B$ are finite and $C$ is malnormal in $A$, i.e.
$C^a\cap C=1$ for every $a\in A\setminus C$. Note that the normalizer condition for $G$ is equivalent to the condition that any finitely generated subgroup of $G$ has finite index in its normalizer.
Moreover, the normalizer condition for a finite graph of finite groups
can be verified algorithmically (see~\cite{B3}).

\begin{thm}~\label{main1} Let $(\mathbb{G},\Gamma)$ be a finite tree of finite groups, which satisfies the normalizer condition. Then its fundamental group $\pi_1(\mathbb{G},\Gamma)$ is subgroup conjugacy separable.
\end{thm}


We deduce these theorems from Theorems~\ref{main3} and~\ref{main2}, and Proposition~\ref{eq}, where the following variation of Definition~\ref{def2} is used.

\begin{defn}\label{SICS}  {\rm
1) For two subgroups $A$ and $B$ of a group $G$, we say that $A$ is {\it conjugate into}~$B$, if there is an element $g\in G$ such that $A^g$ is a subgroup of $B$.

2) A group $G$ is called
{\it subgroup into-conjugacy separable} (abbreviated as SICS) if for any two finitely generated subgroups $H_1,H_2\leqslant G$ such that $H_2$ is not conju\-gate into $H_1$, there exists a homomorphism $\phi$ from $G$ to a finite group $\overline{G}$ such that $\phi(H_2)$ is not conjugate into $\phi(H_1)$ in $\overline{G}$.}
\end{defn}

\begin{prop}\label{eq} Let $G$ be a virtually free group.
Suppose that $G$ is subgroup into-conjugacy separable. Then $G$ is subgroup conjugacy separable.
\end{prop}

\begin{thm}~\label{main3} Free groups are subgroup into-conjugacy separable.
\end{thm}

\begin{thm}~\label{main2} Let $(\mathbb{G},\Gamma)$ be a finite tree of finite groups, which satisfies the normalizer condition. Then its fundamental group $\pi_1(\mathbb{G},\Gamma)$ is subgroup into-conjugacy separable.
\end{thm}


{\bf Methods.} In the proof of Theorem~\ref{main3} we use coverings of graphs, while in the proof of Theorem~\ref{main2}, we use a 3-dimensional topological realization of the graph of groups. Any covering of this realization can be obtained by gluing of some elementary spaces, which we call covering pieces.
This technique is similar to that, which was developed by the first author in the papers~\cite{B3} and~\cite{B2} for a classification of groups with the M.~Hall property.

To construct certain coverings (and so certain subgroups), we use a gluing schema, which comes from Theorem~\ref{graph}
of F\"{u}redi, Lazebnik, Seress, Ustimenko, and Woldar
on the existence of $(r,s)$-bipartite graphs without short cycles.

The structure of the paper is the following. In Section~2 we prove Proposition~\ref{eq}, in Section~3 we give
some auxiliary statements. In Sections~4 and 5 we prove our main Theorems~\ref{main3} and~\ref{main2}.

\section{The SICS-property implies the SCS-property\\ for virtually free groups}~\label{CS}
\begin{lem}~\label{Tak} Let $H_1,H_2$ be two finitely generated subgroups of a virtually free group $G$, such that
$H_1$ is conjugate into $H_2$ and $H_2$ is conjugate into $H_1$.
Then $H_1$ is conjugate to~$H_2$.
\end{lem}

\demo  It is sufficient to prove this theorem for finitely generated $G$.

Let $H_1^{g_1}\leqslant H_2$ and $H_2^{g_2}\leqslant H_1$ for some $g_1,g_2\in G$.
Then $H_1^g\leqslant H_1$ for $g=g_1g_2$.
Moreover, $H_1^g=H_1$ if and only if  $H_1^{g_1}=H_2$ and $H_2^{g_2}=H_1$.

Suppose that $H_1^g$ strictly less than $H_1$.
Then, for $x=g^{-1}$, we have the strictly ascending chain of subgroups: $H_1<H_1^x<H_1^{x^2}<\dots$.
Let $F$ be a free normal subgroup of finite index in $G$.
We compare this chain with the chain $H_1\cap F\leqslant (H_1\cap F)^{x}\leqslant (H_1\cap F)^{x^2}\leqslant \dots$.

Since
the indices $|H_1^{x^i}:(H_1\cap F)^{x^i}|$ are finite and independent of $i$, and since
the indices $|H_1^{x^i}:H_1|$ are increasing with $i$, the second chain is also strictly ascending:
$H_1\cap F< (H_1\cap F)^{x}< (H_1\cap F)^{x^2}< \dots$.
This contradicts to the theorem of M.~Takahasi (see~\cite{Tak}, or~\cite[Theorem  14.1]{KM}),
which claims, that a free group of a finite rank (in our case $F$) does not contain a strictly ascending chain
of subgroups of a finite bounded rank.
Thus, $H_1^g=H_1$ and so $H_1^{g_1}=H_2$.
$\Box$


\medskip

{\it Proof of Proposition~\ref{eq}.} Let $H_1,H_2$ be two non-conjugate, finitely generated subgroups of $G$. By Lemma~\ref{Tak}, w.l.o.g. we may assume that $H_1$ is not conjugate into $H_2$. Since $G$ is a SICS-group, there
exists a homomorphism $\varphi$ from $G$ to a finite group $\overline{G}$, such that $\varphi(H_1)$ is not conjugate into $\varphi(H_2)$ in $\overline{G}$. In particular, $\varphi(H_1)$ is not conjugate to $\varphi(H_2)$ in $\overline{G}$. So, $G$ is subgroup conjugacy separable.
$\Box$


\newpage

\section{Auxiliary statements}

\begin{lem}\label{prop} Let $H_1,H_2$ be subgroups of a group $G$. Then the following conditions are equivalent:

{\rm (1)} $H_2$ is conjugate into every finite index subgroup of $G$, containing $H_1$;

{\rm (2)} For every finite quotient of $G$, the image of $H_2$ is conjugate into the image of $H_1$.
\end{lem}

{\it Proof.} $(2)\Rightarrow (1)$: Suppose that (2) holds and let $D$ be a finite index subgroup of $G$,
containing $H_1$.
Then $D$ contains a finite index subgroup $N$, which is normal in $G$.
By (2), the image of $H_2$ in $G/N$ is conjugate into the image of $H_1$ in $G/N$.
This implies that $H_2$ is conjugate into $H_1N$, and so into $D$.

$(1)\Rightarrow (2)$: Suppose that (1) holds and let $G/N$ be a finite quotient of $G$. By (1), $H_2$ is
conjugate into $H_1N$. Then the image of $H_2$ in $G/N$ is conjugate into the image of $H_1$ in $G/N$.
$\Box$

\begin{lem}\label{Serre3} Let $G$ be a free product: $G=G_1\ast G_2\ast \dots \ast G_l\ast F$, and let $H=\langle h_1,\dots,h_r\rangle$ be a finitely generated subgroup of $G$. Suppose that each $h_i$ and each $h_ih_j$ is conjugate into a factor $G_k$, where $k$ depends on $i$ (on i,j). Then the whole $H$ is conjugate into some $G_s$.\end{lem}

\demo We may assume that $H\neq 1$.
By the Bass-Serre theory (see~\cite{Serre}), $G$ acts on a simplicial tree $T$ without inversions of edges so that the stabilizers of vertices of $T$ are conjugate to $G_1,\dots, G_l,F$. So, each $h_i$ and each $h_ih_j$
stabilize a vertex of $T$. By Corollary~3 in~\cite[Chapter~I, Section~6.5]{Serre}
of Serre, $H$ stabilizes a vertex of $T$ and hence $H$ is conjugate into some $G_s$ or into $F$.
The last cannot happen, since $H$ contains a non-trivial generator $h_i$, which is conjugate into some $G_k$.  $\Box$

\medskip

A graph $K$ is called {\it bipartite}
if the set of its vertices is a disjoint union of two nonempty sets $V_1$ and $V_2$, such that every edge of $K$ connects
a vertex from $V_1$ to a vertex from $V_2$. A bipartite graph is said to be {\it bi-regular} if there exist integers $r,s$ such that
${\text{\rm deg}}(x)=r$ for all $x\in V_1$ and ${\text{\rm deg}}(y)=s$ for all $y\in V_2$.
In this case $(r,s)$ is called {\it bi-degree} of $K$. Note that the lengths of cycles in a bipartite graph are always even.

\begin{thm}{\rm \cite{FL}.}\label{graph} For any natural $r,s,t\geqslant 2$, there exists a finite connected bipartite graph of bi-degree $(r,s)$, with length of smallest cycle exactly $2t$.
\end{thm}
An $r$-{\it star} is a tree with $r+1$ vertices and $r$ edges, outgoing from one common vertex.
This vertex will be called {\it central} and the other ones {\it peripherical}.
It is convenient to reformulate a weaker version of this theorem.

\begin{thm}\label{refgraph} For any natural $r,s,t\geqslant 1$, one can glue several $r$-stars to several $s$-stars, so
that all peripherical vertices of $r$-stars will be identified (by some bijection) with all peripherical vertices of $s$-stars, the resulting graph will be connected, and
it will not have cycles of length smaller than $t$.
\end{thm}

\section{SICS-property for free groups}

\subsection{Notations}

Our proof of Theorem~\ref{main0} uses coverings of labeled graphs.
Here we define a core of a covering, an outer edge, and an outer vertex of a core.

Let $\Gamma$ be a graph. By $\Gamma^0$ we denote the set of its vertices and by $\Gamma^1$ the set of its edges.
The inverse of an edge $e\in \Gamma^1$ is denoted by $\overline{e}$, the initial and the terminal vertices of $e$ are denoted by $i(e)$ and $t(e)$.

Let $F$ be a free group with finite basis $x_1,\dots,x_n$.
Let $R$ be the graph consisting of one vertex $v$ and $n$ oriented edges $e_1,\dots,e_n$.
We label $e_i$ by $x_i$ and $\overline{e}_i$ by $x_i^{-1}$.
We will identify $F$ with $\pi_1(R,v)$ by identifying $x_i$ with the homotopy class $[e_i]$.

To every subgroup $H\leqslant F$ corresponds a covering map $\varphi: (\Gamma_H,v_H)\rightarrow (R,v)$,
such that $H$ is the image of the induced map $\varphi_{\ast}:\pi_1(\Gamma_H,v_H)\rightarrow \pi_1(R,v)$.
We lift the labeling of $R$ to $\Gamma_H$. So, an edge $e$ of $\Gamma_H$ is labeled  by $x$ if its image $\varphi(e)$ is labeled by $x$.

If $H$ is finitely generated, then $\Gamma_H$ has a finite {\it core}, ${\text{\rm Core}}(\Gamma_H)$, that is a finite connected subgraph, which is homotopy equivalent to $\Gamma_H$. We can enlarge ${\text{\rm Core}}(\Gamma_H)$ if necessary and assume that $v_H$ is a vertex of ${\text{\rm Core}}(\Gamma_H)$ and that every vertex of ${\text{\rm Core}}(\Gamma_H)$ has valency 1 or $2n$. The vertices of valency 1 and the edges incident to these vertices
are called {\it outer}.
All other vertices and edges of ${\text{\rm Core}}(\Gamma_H)$ are called {\it inner}.
Let $e$ be an oriented outer edge of ${\text{\rm Core}}(\Gamma_H)$, which starts at an outer vertex. Then there is a unique oriented path $e_1e_2\dots e_k$
in ${\text{\rm Core}}(\Gamma_H)$, such that $e_1=e$, the labels of edges $e_i$ are coincide, and the last edge $e_k$ is outer.
We will write $e_1^{\,\,\ast}=e_k$. Clearly $\overline{e}_k^{\,\,\ast}=\overline{e}_1$.
Thus, we get a free involution ${\ast}$ on the set of outer edges of ${\text{\rm Core}}(\Gamma_H)$.


\subsection{Proof of Theorem~\ref{main3}}

Let $F$ be a free group with finite basis $x_1,\dots,x_n$.
Let $H_1,H_2$ be two nonconjugate finitely generated subgroups of $F$ such that $H_2$ is not conjugate into $H_1$.
By Lemma~\ref{prop}, it is sufficient to construct a finite index subgroup $D$ of $F$, that contains $H_1$ and does not contain a conjugate of $H_2$.

Let $H_2=\langle h_1,\dots, h_r\rangle$ and let $C=2 \max\{ |h_1|,\dots ,|h_r|\}$. Since $F$ is a residually finite,
there exists a normal subgroup $K$ of finite index in $F$, such that $K$ does not contain nontrivial elements of $F$ of length $C$ or smaller. Since $K$ is normal, $K$ does not contain any conjugate to these elements.
This means that the covering graph $\Gamma_{K}$ is finite, every its vertex has valency $2n$, and

$${\text {\it every cycle in}}\hspace*{2mm} \Gamma_{K}\hspace*{2mm} {\text {\it has length at least}}\hspace*{2mm} C+1.\eqno{(1)}
$$

\medskip

Without loss of generality, we assume that the vertices of ${\text{\rm Core}}(\Gamma_{H_1})$ have valency 1 or $2n$.
Now we will embed ${\text{\rm Core}}(\Gamma_{H_1})$ into a finite labeled graph $\Delta$ without outer edges.
Let ${\mathcal{E}}$ be the set of all edges of ${\text{\rm Core}}(\Gamma_{H_1})$, that start at outer vertices of ${\text{\rm Core}}(\Gamma_{H_1})$.
For every edge $e\in \mathcal{E}$ we choose an edge $\widehat{e}$ in $\Gamma_{K}$ with the same label.
Let $\Delta$ be the labeled graph, obtained from the disjoint union of graphs


$${\text{\rm Core}}(\Gamma_{H_1})\bigsqcup\,\, \underset{e\in \mathcal{E}}{\sqcup}\,\, \Bigl(\Gamma_{K}\setminus\{\widehat{e},\,\overline{\widehat{e}}\}\Bigr)\eqno{(2)} $$
by identifying the vertices $\alpha(e)$ with $\alpha(\widehat{e})$ and $\omega(e^{\ast})$
with $\omega(\widehat{e})$ for every $e\in \mathcal{E}$.

\hspace*{35mm}\includegraphics[scale=0.4]{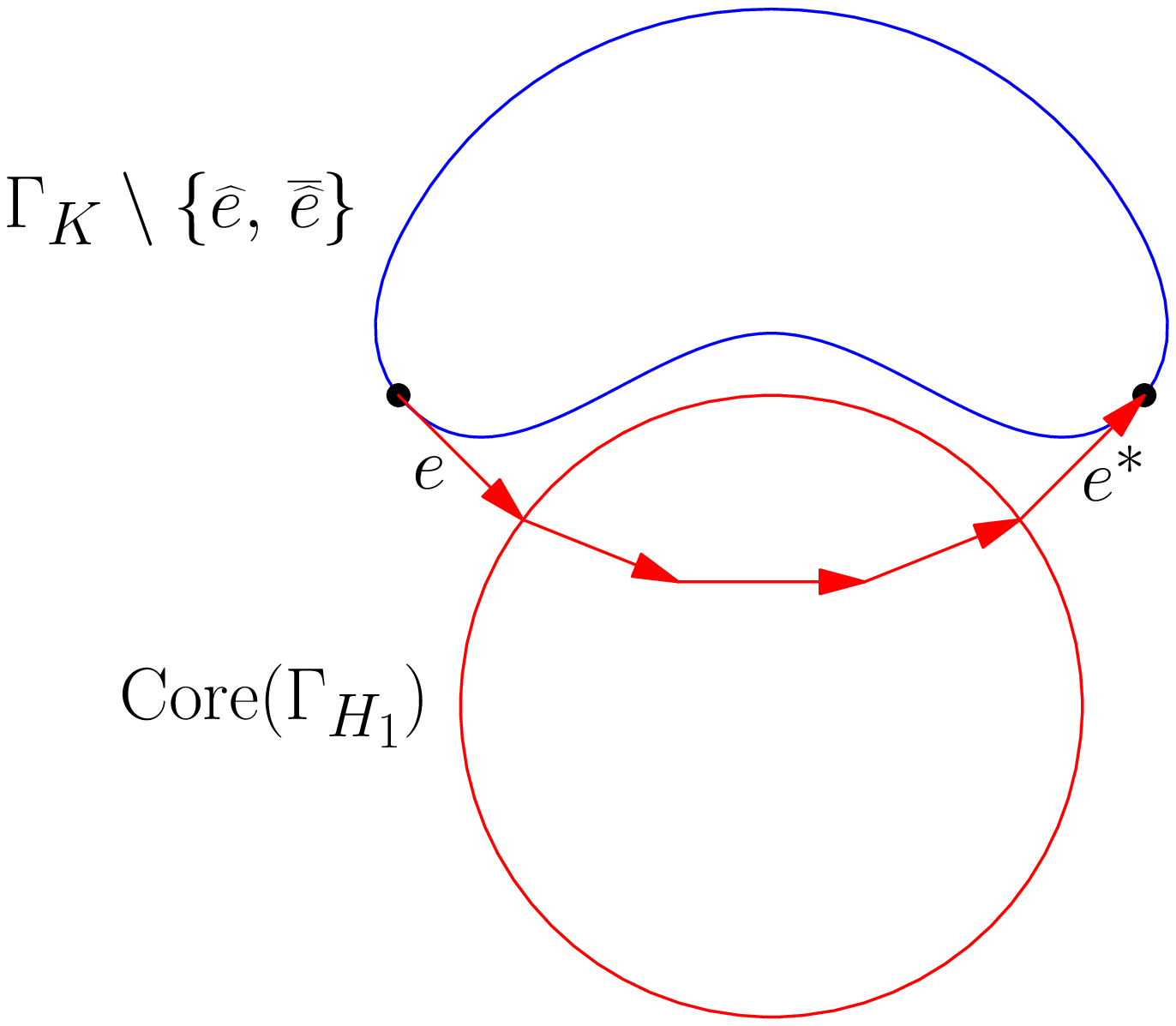}

\vspace*{-66mm}
\begin{center} Figure 1
\end{center}

Since every vertex of $\Delta$ has valency $2n$,
there is a finitely sheeted covering map $\psi: (\Delta, v_{H_1})\rightarrow (R,v)$, respecting the labeling.
Thus $\Delta=\Gamma_D$ for some finite index subgroup $D$ of $F$. Since $\Gamma_{H_1}$ is a subgraph of $\Gamma_D$, the subgroup $D$ contains $H_1$ as a free factor: $D=H_1\ast L$.

We show, that $H_2=\langle h_1,\dots ,h_r\rangle$ is not conjugate into $D$.
Suppose the contrary: $H_2^g\leqslant D$ for some $g\in F$, Then
every element $h\in \{h_i^g$, $(h_ih_j)^g: i,j=1,\dots, r\}$ can be represented by a closed path $l(h)$ in
$\Delta$ based at $v_{H_1}$. By definition of the constant $C$, every path $l(h)$
can be freely homotopic to a closed path in $\Delta$ of length at most~$C$.
By Condition (1) and by Construction (2), every such path is freely homotopic to a closed path in ${\text{\rm Core}}(\Gamma_{H_1})$.
This means that every element $h\in \{h_i^g$, $(h_ih_j)^g: i,j=1,\dots, r\}$ can be conjugated into
$H_1$ by an element $d(h)\in D$.
By Lemma~\ref{Serre3}, $H_2^g$ can be conjugated into $H_1$ by an element $d\in D$.
This contradicts to the assumption, that $H_2$ cannot be conjugated into $H_1$ in $F$. $\Box$

\section{SICS-property for virtually free groups}

By~\cite{BKS},
every finitely generated virtually free group is the fundamental group of a finite graph of finite groups (see also~\cite[Chapter IV, Theorem 1.6]{DD} and historical comments on page 133 of~\cite{DD}). We will also represent these groups as fundamental groups of some graphs of spaces (3-dimensional complexes). Below we introduce notations and recall some definitions. In Subsection~\ref{mainproof} we prove Theorem~\ref{main2}.

\subsection {Graphs of groups}\label{1}


A graph of groups $(\mathbb{G},\Gamma)$ is a system consisting of a connected graph $\Gamma$, of {\it vertex groups} $G_v$, $v\in \Gamma^0$, of {\it edge groups} $G_e$, $e\in \Gamma^1$, and of {\it boundary monomorphisms} $\rho_e^{i}:G_e\rightarrow G_{i(e)}$ and $\rho_e^{t}:G_e\rightarrow G_{t(e)}$, $e\in \Gamma^1$, which satisfy
$G_e=G_{\overline{e}}$ and
$\rho_{e}^{i}=\rho_{\overline{e}}^{t}$.

\medskip



A {\it path} in the graph of groups $(\mathbb{G},\Gamma)$ is a sequence of the form $g_1e_1g_2e_2\dots e_kg_{k+1}$,
where $e_1e_2\dots e_k$ is a path in $\Gamma$, $g_s\in G_{i(e_s)}$ and $g_{s+1}\in G_{t(e_s)}$ for $s=1,2,\dots ,k$.
This path is {\it closed}, if $i(e_1)=t(e_{k+1})$; in this case we say that it is {\it based} at the vertex $i(e_1)$. There is a usual (partial) multiplication of paths in $(\mathbb{G},\Gamma)$.

Now we define three types of {\it elementary transformations} of a path $l=g_1e_1g_2e_2\dots e_kg_{k+1}$:

\medskip

1) replace a subpath of $l$ of the form $aeb$, where $e\in \Gamma^1$, $a\in G_{i(e)}$ and $b\in G_{t(e)}$,
by the path $a_1eb_1$, where $a_1=a(\rho_{e}^{i}(g))^{-1}$ and $b_1=\rho_{e}^{t}(g)b$ for some $g\in G_{i(e)}$;

2) replace a subpath of $l$ of the form $ae1\overline{e}b$, where $e\in \Gamma^1$ and $a,b\in G_{i(e)}$, by the element $ab\in G_{i(e)}$;

3) this is the transformation inverse to 2).

\medskip

Two paths $l$ and $l'$ in $(\mathbb{G},\Gamma)$ are called {\it equivalent}, if $l'$ can be obtained from $l$
by a finite number of elementary transformations. The equivalence class of $l$ is denoted by $[l]$.


{\it The fundamental group of the graph of groups} $(\mathbb{G},\Gamma)$ with respect to a vertex $v\in \Gamma^0$, denoted $\pi_1(\mathbb{G},\Gamma,v)$, is the set of equivalence classes of all closed paths in $(\mathbb{G},\Gamma)$ based at $v$ with respect to the multiplication $[l_1][l_2]=[l_1l_2]$.

Denote $G=\pi_1(\mathbb{G},\Gamma,v)$. Every element $g\in G$ can be represented by a closed path
$g_1e_1g_2e_2\dots e_kg_{k+1}$ with minimal $k=k(g)$. We call such $k$ the {\it length} of $g$ and denote it by $|g|$.

Note, that every vertex group $G_u$ of the graph of groups $(\mathbb{G},\Gamma)$ can be embedded into
$G$ by the following rule. Choose a path $p$ is $\Gamma$ from $v$ to $u$. The map $g\mapsto [pgp^{-1}]$, $g\in G_u$ determines
an embedding of $G_u$ into $G$.  If we choose another path from $v$ to $u$, the resulting subgroup will be conjugate to the first one.
Thus, $G_u$ canonically determines the conjugacy class of a subgroup of $G$. Any subgroup of this class will be called
a {\it vertex subgroup of $G$}, corresponding to $G_u$.

\subsection{Graph of spaces}\label{2}

Below all spaces are assumed to be path connected topological spaces.
In particular, their fundamental groups are well defined (up to isomorphism).

A {\it graph of spaces} $(\mathbb{X},\Gamma)$ is a system consisting of a connected graph $\Gamma$, of {\it vertex spaces} $X_v$, $v\in \Gamma^0$, of {\it edge spaces} $X_e$, $e\in \Gamma^1$, and of $\pi_1$-injective continuous {\it boundary maps} $\partial_e^{i}:X_e\rightarrow X_{i(e)}$ and $\partial_e^{t}:X_e\rightarrow X_{t(e)}$, $e\in \Gamma^1$, which satisfy
$X_e=X_{\overline{e}}$ and
$\partial_e^{i}=\partial_{\overline{e}}^{t}$. For the later it is convenient to think, that $X_e$ and $X_{\overline{e}}$ are two copies of the same space.



The {\it topological realization} of the graph of spaces $(\mathbb{X},\Gamma)$, denoted ${\text{\rm Real}}(\mathbb{X},\Gamma)$,
is defined to be the quotient obtained from
$$\underset{v\in \Gamma^{0}}{\coprod}X_v\,\, {\sqcup} \,\, \underset{e\in \Gamma^1}{\coprod} (X_{e}\times [0,1]),$$
by gluing $(x,0)$ to $\partial_e^{i}(x)$ and $(x,1)$ to $\partial_e^{t}(x)$ for every
$e\in \Gamma^1$ and $x\in X_e$, and by identifying the spaces $X_e\times [0,1]$ and $X_{\overline{e}}\times [0,1]$
through the map $(x,t)\rightarrow (x,t-1)$, $x\in X_{e}$, $t\in [0,1]$.
Denote $X={\text{\rm Real}}(\mathbb{X},\Gamma)$.



\medskip

\vspace*{-0.5cm}

{\hspace*{1.7cm}\includegraphics[scale=0.5]{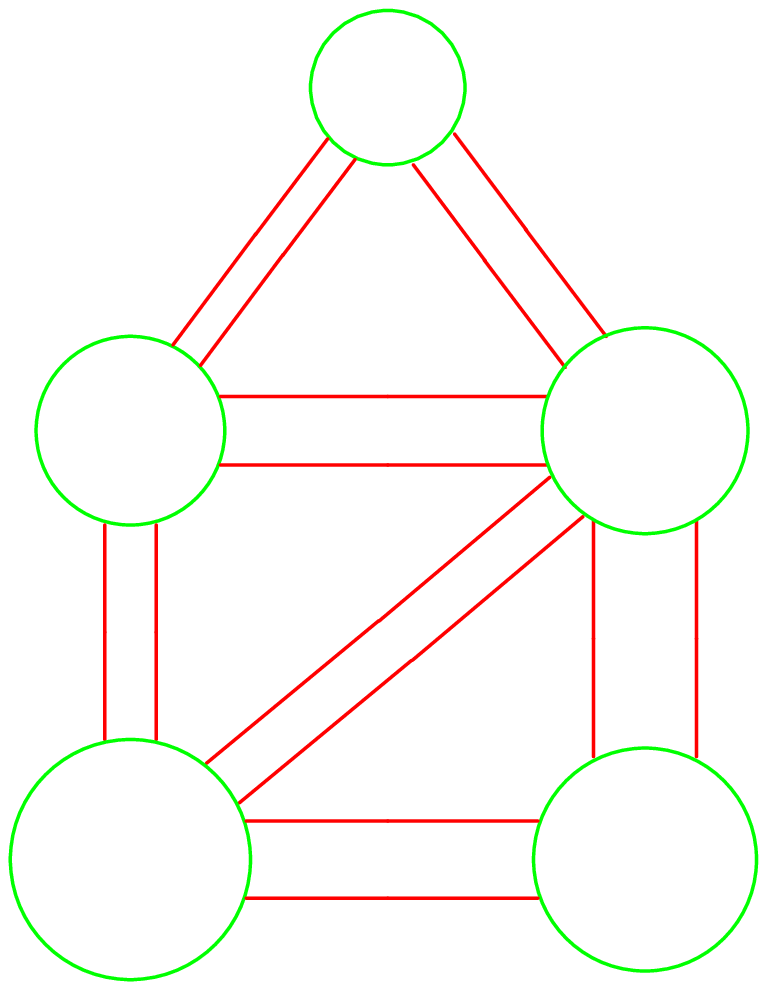}}

\vspace*{-15cm}

{\hspace*{8cm}\includegraphics[scale=0.5]{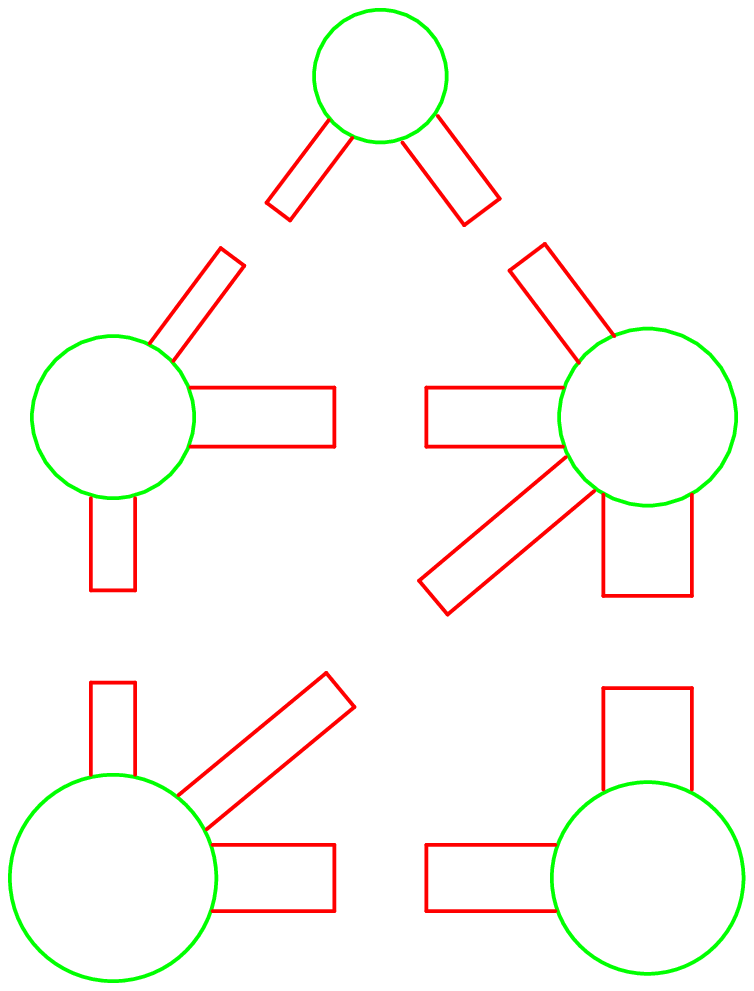}}

\vspace*{-7.9cm}

\hspace*{30mm} {\small The space $X$ \hspace*{37mm} The pieces of $X$}

\vspace*{-4mm}
\begin{center}
Figure 2\hspace*{50mm}Figure 3
\end{center}




A {\it body} of $X$ is a subspace of $X$ of the form $X_v$, $v\in \Gamma^0$.
The {\it piece} associated with the body~$X_v$, denoted $\mathcal{N}(X_v)$, is defined to be the quotient obtained from
the topological space
$$X_v\,\, {\sqcup} \,\, \underset{e\in \Gamma^1, i(e)=v}{\coprod} (X_{e}\times [0, \frac{1}{2}]),$$
by identifying $(x,0)$ with $\partial_e^{i}(x)$ for each edge $e\in \Gamma^1$ outgoing from $v$
and for every $x\in X_e$. The subspaces $X_{e}\times [0, \frac{1}{2}]$
are called the {\it handles} of this piece. The subspaces $X_e\times \{\frac{1}{2}\}$ are called the {\it faces} of this piece.


Any covering of a piece in $X$ is called a {\it covering piece}. Lifts of the body, of the handles, and of the faces of the piece are called the body, the handles, and the faces of the covering piece.
Note that if $\Gamma$ is finite and $\pi_1(X_v)$ is finite for every $v\in \Gamma^0$,
then there is only a finite number of covering pieces, up to homeomorphism.

Clearly, the space $X$ can be obtained by an appropriate gluing of all pieces $\mathcal{N}(X_v)$, $v\in \Gamma^{0}$, along their free faces.
Every covering space of $X$ can be obtained by gluing of (may be infinitely many) copies of covering pieces along their faces.

A topological space $Z$ is called a {\it pre-covering} of $X$,
if $Z$ is a connected subspace of some covering of $X$,
which can be presented as a result of gluing of some covering pieces along their faces.
A face $F\subseteq Z$ is called a {\it free face} of $Z$, if it is a face of exactly one handle of~$Z$.
A handle of $Z$ which contains a free face is called a {\it free handle} of $Z$.

With every pre-covering $Z$ of $X$ we can naturally associate a graph $\Delta$
by collapsing its bodies to vertices and its handles to ``half''-edges.
Let $p:Z\rightarrow \Delta$ be the collapsing map for $Z$.
We equip $Z$ with the pseudometric induced by the usual path metric on~$\Delta$
(where the ``half''-edges have length $\frac{1}{2}$).
In particular, the distance between any two points of a body of $Z$ is zero
and the maximal distance between two points of a handle is $\frac{1}{2}$.

Let $\Delta^{\frac{1}{2}}$ be the set of middle points of edges of length 1 in $\Delta$.
A curve $c:[0,1]\rightarrow \Delta$ is called {\it regular}, if $c$ has endpoints in
$\Delta^0\cup \Delta^{\frac{1}{2}}$ and is locally injective on
$[0,1]\setminus c^{-1}(\Delta^0\cup \Delta^{\frac{1}{2}})$.
The $e$-{\it length} of a regular curve $c$ in $\Delta$, denoted $|c|_e$, is the sum of lengths of the ``half''-edges which $c$ passes.

A curve $\gamma: [0,1]\rightarrow Z$ is called {\it regular}, if the curve $p\circ \gamma:[0,1]\rightarrow \Delta$
is regular. The $e$-{\it length} of a regular curve $\gamma$ in $Z$, denoted $|\gamma|_e$, is defined to be
$|p\circ \gamma|_e$.
Roughly speaking, $|\gamma|_e$ is the number of handles which $\gamma$ passes, divided by 2.






\subsection{From graphs of groups to graphs of spaces}\label{3}

With every group $G$ we associate the 2-dimensional CW-complex $Space(G)$, consisting of the unique vertex ${\bold{u}}_G$, the edge set $\{e_g\,|\, g\in G\}$, and the set of 2-cells $\{D_{a,b}\,|\, a,b\in G\}$, where the boundary of $D_{a,b}$ is glued along the path $e_ae_b\overline{e_{ab}}$.
We identify the groups $G$ and $\pi_1(Space(G),{\bold{u}}_G)$ through the canonical isomorphism $g\mapsto [e_g]$.

With every embedding of groups $\varphi:H_1\hookrightarrow H$ we associate the embedding of complexes
$Space(H_1)\hookrightarrow Space(H)$, such that the induced homomorphism of fundamental groups coincides with $\varphi$.

Now, with any graph of groups $(\mathbb{G},\Gamma)$ we associate the graph of spaces
$(\mathbb{X},\Gamma)$,
such that $X_w=Space(G_w)$ for $w\in \Gamma^0\cup \Gamma^1$ and
the embeddings of spaces, $\partial_e^i$, $\partial_e^t$, correspond to the embeddings of groups $\rho_e^i$, $\rho_e^t$.
For any vertex $v\in \Gamma$,
the groups $\pi_1(\mathbb{G},\Gamma,v)$  and $\pi_1({\text{\rm Real}}(\mathbb{X},\Gamma),{\bold{u}}_{G_v})$
are canonically isomorphic and we will identify them through this isomorphism.
Every element $g\in \pi_1(\mathbb{G},\Gamma,v)$ can be realized by a
regular closed path $\gamma(g)$ in ${\text{\rm Real}}(\mathbb{X},\Gamma)$ based at ${\bold{u}}_{G_v}$, such that the number of subspaces $X_e\times [0,1]$ it crosses is equal to $|g|$; in our notations we have $|\gamma(g)|_e=|g|$.


\subsection{Subgroups and coverings, cores of coverings}

Simplifying notations in the previous section, we denote $G=\pi_1(\mathbb{G},\Gamma,v)$,  $X={\text{\rm Real}}(\mathbb{X},\Gamma)$, and $x={\bold{u}}_{G_v}$.  We may assume $G=\pi_1(X,x)$.

For every subgroup $H$ of $G$ there exists a covering map $\psi:(Y,y)\rightarrow (X,x)$, such that $H=\psi_{\ast}(\pi_1(Y,y))$.
The space $Y$ can be presented as the topological realization of a graph of spaces $(\mathbb{Y},\Delta)$. The vertex spaces and the edge spaces of $(\mathbb{Y},\Delta)$ are connected components of $\psi^{-1}(X_u)$, $u\in \Gamma^0$, and of $\psi^{-1}(X_e)$, $e\in \Gamma^1$ respectively.

If $H$ is finitely generated, then there is a pre-covering $Y_0\subseteq Y$ of $X$, such that\\
1) $Y_0$ can be obtained by gluing of finitely many covering pieces along some of their faces;\\
2) every loop in $Y$ can be freely homotoped into $Y_0$.

\medskip

Any such pre-covering will be called a {\it core} of $Y$. We choose one of them and denote it by $Core(Y)$.
The closure of every connected component of  $Y\setminus Core(Y)$ will be called a {\it thick tree}.
Every thick tree grows from a free face of the core and it can be deformationally retracted onto this face.

\subsection{Trivial handles}

A covering handle will be called {\it trivial} if its fundamental group is trivial.
The following is a preparation to the proof of Theorem~\ref{main2}.

\subsubsection{A linear order on the set of trivial covering handles}

From now on we assume that $\Gamma$ is a tree.
For every edge $e\in \Gamma^1$, let ${\text{\rm Comp}}(e)$ denote the connected component of $\Gamma\setminus \{e,\overline{e}\}$ which contains $t(e)$.
 For every $m\geqslant 0$ we set $$\Gamma^1(m)=\{e\in \Gamma^1\mid {\text{\rm the number of geometric edges in}}\hspace*{2mm} {\text{\rm Comp}}(e)\hspace*{2mm}{\text{\rm is equal to}}\hspace*{2mm} m\}.$$

Now we choose an arbitrary linear order on each $\Gamma^1(m)$ and extend these orders to a linear
order $\prec$ on $\Gamma^1$ by saying that edges in $\Gamma^1(m_1)$ are smaller than edges in $\Gamma^1(m_2)$
if and only if $m_1<m_2$.

We will consider covering spaces up to equivalence of coverings.
Let $\mathcal{C}_1, \mathcal{C}_2$ be two trivial covering handles, which cover the handles $X_{e_1}\times [0,\frac{1}{2}]$ and $X_{e_2}\times [0,\frac{1}{2}]$ of $X$.
We say that $\mathcal{C}_1$ is {\it smaller} than $\mathcal{C}_2$ and write $\mathcal{C}_1\prec_{\ast} \mathcal{C}_2$, if $e_1\prec e_2$.





\subsubsection{Extensions of pre-coverings through their free trivial handles}

We explain a construction, which enables to extend pre-coverings of $X$ through their trivial free handles keeping the fundamental group unchanged.

Let $Y$ be a pre-covering of $X$ and let $A$ be a trivial free handle of $Y$.
 By definition, the face of $A$ is free in $Y$ and $A$ is the universal covering of a handle in $X$, say $C_1$.
Let $C_2$ be the other handle in $X$ which has a common face with $C_1$, and let $\mathcal{N}$ be the piece in $X$ which contains $C_2$.
Consider the universal covering $\widetilde{\mathcal{N}}\rightarrow \mathcal{N}$ and some lift $B$ of $C_2$ in $\widetilde{\mathcal{N}}$.
We can extend $Y$ by gluing $Y$ and  $\widetilde{\mathcal{N}}$ along the free faces of $A$ and $B$.
So we obtain a new pre-covering of $X$ with the same fundamental group as $Y$.

{\hspace*{27mm}\includegraphics[scale=0.5]{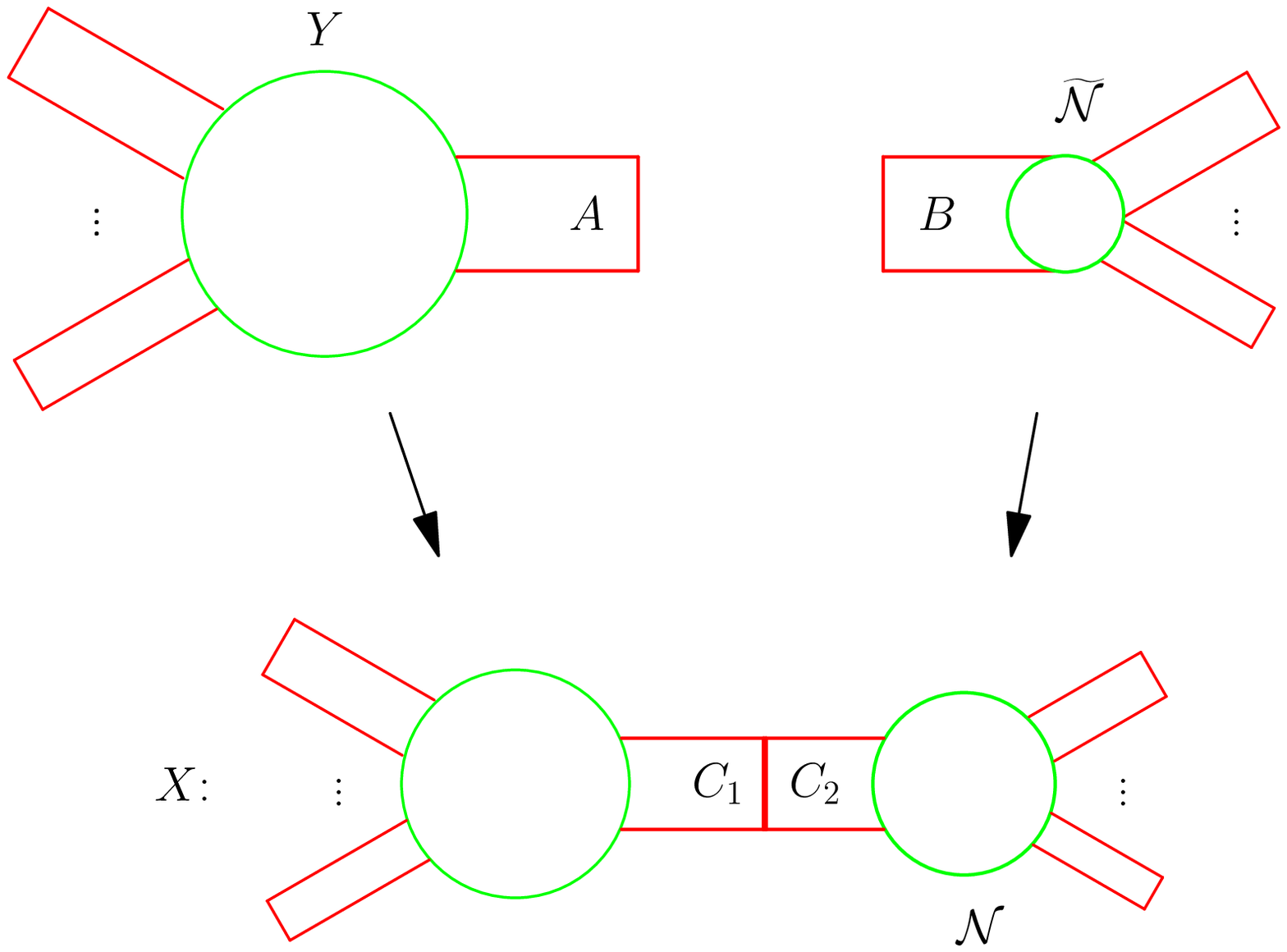}}

\vspace*{-8.4cm}

\begin{center} Figure 4
\end{center}

As a preparation to the proof of Theorem~\ref{main2}, we explain a more complicated gluing.
Let $B_1,\dots,B_s$ be all lifts of $C_2$ in $\widetilde{\mathcal{N}}$. Now we take $s$ copies
of $Y$, say $Y_1,\dots ,Y_s$ (with the handles $A_1,\dots,A_s$ corresponding to $A$) and glue them to
$\widetilde{\mathcal{N}}$ so that the free face of $A_i$ will be glued to the free face of $B_i$ for $i=1,\dots,s$.
The resulting space is again a pre-covering of $X$ and its fundamental group
is isomorphic to the free product of $s$ copies of $\pi_1(Y)$.
We call the tuple $(B_1,\dots,B_s; \widetilde{\mathcal{N}}$) the {\it tuple associated with the handle~$A$}.

The following lemma allows to prove~Theorem~\ref{main2} by induction.

\begin{lem}\label{order} Let $K$ be a handle of $\widetilde{\mathcal{N}}$, different from $B_1,\dots ,B_s$. Then
$K\prec_{\ast} A$.
\end{lem}

{\it Proof.}
The proof is straightforward and uses the assumption, that $\Gamma$ is a tree. $\Box$

\subsection{Proof of Theorem~\ref{main2}}\label{mainproof}

Let $G$ be the fundamental group of a finite graph of finite groups: $G=\pi_1(\mathbb{G},\Gamma,v)$,
where $\Gamma$ is a tree and suppose that $G$ satisfies  the normalizer condition.
By Section~\ref{3}, we can write
$G=\pi_1(X,x)$, where
$X$ is the topological realization of the graph of spaces $(\mathbb{X},\Gamma)$ associated with the graph of groups $(\mathbb{G},\Gamma)$.



Let $H_1,H_2$ be two finitely generated subgroups of $G$ and suppose that
$H_2$ is not conjugate into $H_1$ in $G$. We will construct a finite index subgroup $H_3$ of $G$,
such that $H_1$ is contained in $H_3$ and $H_2$ is not conjugate into $H_3$. Then Theorem~\ref{main2} will
immediately follow from Lemma~\ref{prop}.

Suppose that $H_2=\langle h_1, h_2, \dots, h_k\rangle$ and let $C=2\max\{|h_i|: i=1,\dots ,k\}+3$, where $|\!\cdot\!|$ is the length function on $G$ (see Section~\ref{1}).

Consider the covering $\varphi:(Y,y)\rightarrow (X,x)$ which corresponds to $H_1$, i.e.
$H_1=\varphi_{\ast}(\pi_1(Y,y))$.
By enlarging the core of $Y$ if necessary, we may assume that $y\in Core(Y)$.
We will complete $Core(Y)$ to a finite sheeted covering space $Z$ and put then $H_3=\pi_1(Z,y)$.
We will do that in several steps.

1) Since $G$ satisfies the normalizer condition, there is a compact pre-covering $Y_0\subseteq Y$ that contains $Core(Y)$ and whose free faces are trivial.

\medskip

2) Consider the compact pre-covering $Y_1\subseteq Y$ that contains $Y_0$ and whose free faces are at distance $C$ from $Y_0$.
The components of $Y_1\setminus Y_0$ are parts of thick trees, which grow from the free faces of $Y_0$.
By 1) these componens are contractible.

{\hfill \includegraphics[width=12.6cm]{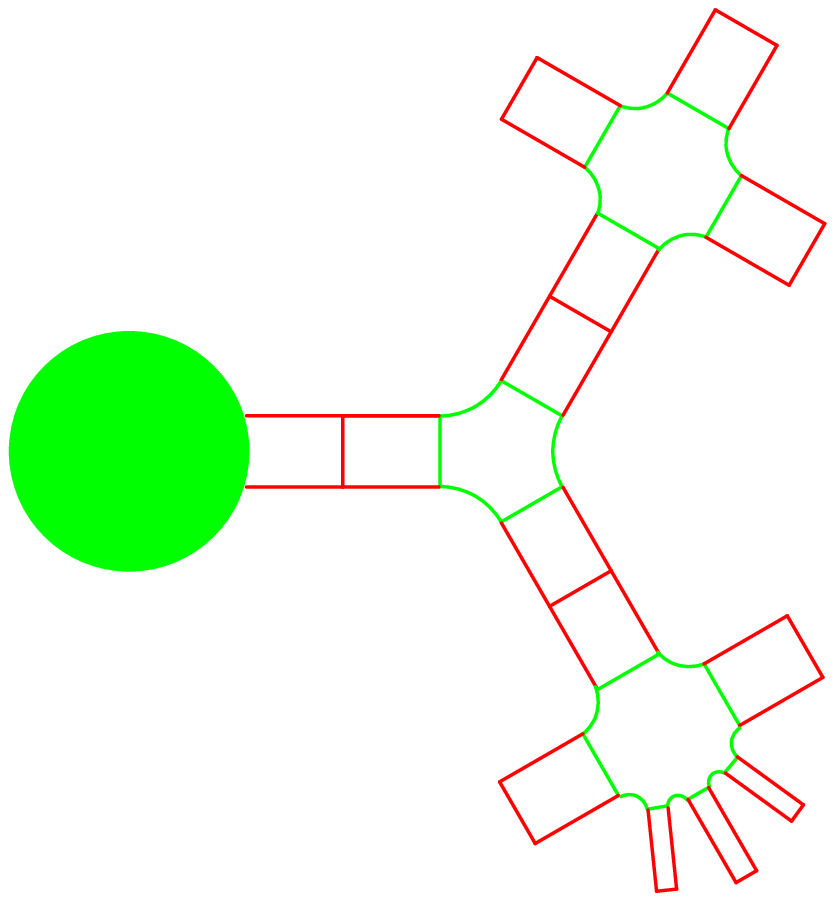}}

\vspace*{-11.5cm}

\begin{center} Figure 5
\end{center}

Note for step 4), that $H_1=\varphi_{\ast}(\pi_1(Y_1,y))$.

\medskip

3) Let $A_1,\dots ,A_r$ be the free handles of $Y_1$ that are maximal, with respect to $\prec_{\ast}$, among all free handles of $Y_1$.
Let $(B_1,\dots,B_s; \widetilde{\mathcal{N}})$ be the tuple associated with the handle $A_1$ (and so with each $A_i$).
Taking into account only these handles, it is convenient to think that $Y_1$ has the form of the $r$-star and $\widetilde{\mathcal{N}}$ has the form of the $s$-star.

We take several copies of $Y_1$ and several copies of $\widetilde{\mathcal{N}}$ and
glue them according to Theorem~\ref{refgraph}, where we put $t=C$.

\includegraphics[width=13cm]{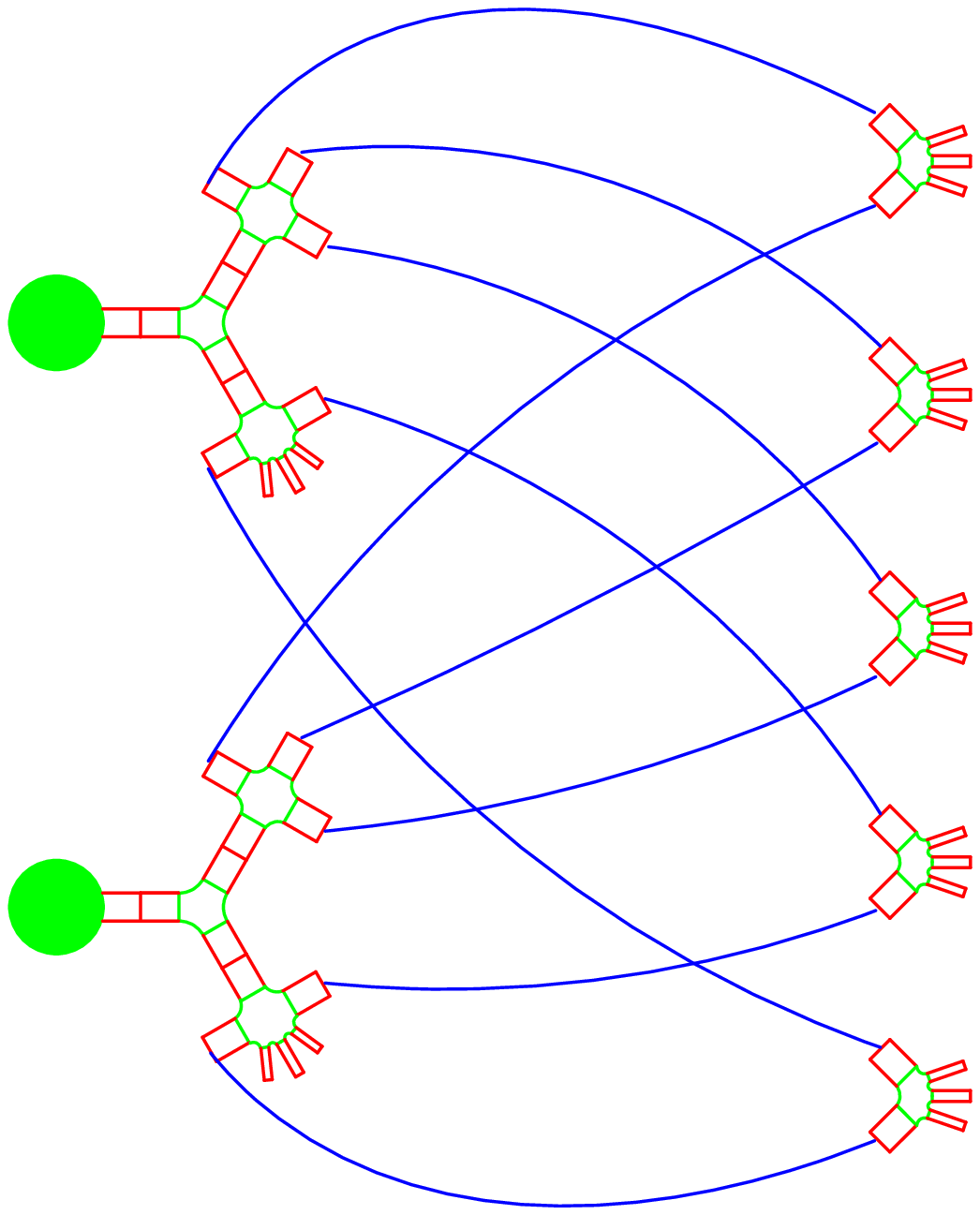}

\vspace*{-9.5cm}

\begin{center} Figure 6
\end{center}


The resulting space $Y_2$ has the following properties.

\medskip

(a) The free handles of $Y_2$ are trivial.

(b) The maximal (with respect to $\prec_{\ast}$) free handles of $Y_2$ are smaller than that of $Y_1$.

\medskip

Property (a) follows from 1).
Property (b) follows from the fact, that the maximal free handles of the copies of $Y_1$ became non-free after gluing them with the free handles of the copies of $\widetilde{\mathcal{N}}$.
Moreover, the handles of copies of $\widetilde{\mathcal{N}}$, which remain free in $Y_2$ are smaller than $A_1$ with respect to $\prec_{\ast}$ by Lemma~\ref{order}.

\medskip

4) We construct $Y_3$ by fulfilling Step 3) for $Y_2$ instead of $Y_1$.
Continuing in this way, we obtain the sequence of regular spaces $Y_1,Y_2,Y_3,\dots $, with the property that
the maximal free handles of $Y_{i+1}$ are smaller than those of $Y_i$. Since the order $\prec_{\ast}$ is finite,
the sequence $Y_1,Y_2,Y_3,\dots $ is finite and the last space, denote it by $Z$, has no free handles.
Then $Z$ is a finite sheeted covering of $X$. Let $H_3$ be the corresponding finite index subgroup of $G$.

Note, that $Z$ is the result of gluing of several copies of $Y_1$, say $Y_{1,1},\dots, Y_{1,n}$, and several
contractible spaces, say $\widetilde{\mathcal{N}}_1,\dots, \widetilde{\mathcal{N}}_m$, along their free faces. So, we have
$\pi_1(Z)\cong\pi_1(Y_{1,1})\ast \dots \ast \pi_1(Y_{1,n})\ast F$
for some free group $F$.
This means, that

$$H_3=H_1^{g_1}\ast \dots \ast H_1^{g_n}\ast F\eqno{(3)}$$
for some $g_1,\dots ,g_n\in G$. After renumbering, we may assume that $Y_1=Y_{1,1}$, and so $g_1=1$.


\medskip

\begin{lem} Every regular loop $\gamma$ in $Z$, which crosses a face of some $\widetilde{\mathcal{N}}_i$
and has the e-length smaller than $C$ is contractible.
\end{lem}

{\it Proof.}
Let $l$ be the minimal natural number, such that $\gamma$ lies in a copy of $Y_l$,
which was used in the construction of $Z$.
For simplicity, we assume that this copy coincides with~$Y_l$.
Note that the $e$-length of any closed regular curve in $Z$
is a nonnegative {\it integer} number (see Section~\ref{2}).

We will proceed by induction on $(l, |\gamma|_e)$,
assuming that the pairs are lexicographically ordered. If $|\gamma|_e=0$, then $\gamma$
entirely lies in a face of $Z$, and hence is contractible. So, we assume that $|\gamma|_e>0$.

{\it Case 1.} Suppose that $l=1$. Then $\gamma$ lies in $Y_1$ and crosses a free face of $Y_1$. Recall that by Step 2),
the distance between any free face of $Y_1$ and $Y_0$ is $C$ and
that every component of $Y_1\setminus Y_0$ is contractible.
Since $|\gamma|_e<C$, the curve $\gamma$ lies in some component of $Y_1\setminus Y_0$ and so is contractible.


\medskip

{\it Case 2.} Suppose that $l\geqslant 2$. Recall that $Y_l$ is the result of gluing of several copies of $Y_{l-1}$
and several $\widetilde{\mathcal{N}}$-spaces. We will say that these $\widetilde{\mathcal{N}}$-spaces have {\it level} $l$.

Since $l$ is minimal, $\gamma$ crosses some $\widetilde{\mathcal{N}}$-space of level $l$, say $\widetilde{\mathcal{N}}_1$. If $\gamma$ completely
lies in $\widetilde{\mathcal{N}}_1$, it is contractible. So, we assume that
$\gamma$ crosses a free face of $\widetilde{\mathcal{N}}_1$, say $\mathcal{F}_1$, and enters into a copy of $Y_{l-1}$,
say $Y_{l-1,1}$.
After running inside $Y_{l-1,1}$ it must again cross a free face $\mathcal{F}_2$ of $Y_{l-1,1}$.

{\it Subcase 2.1.} Suppose that $\mathcal{F}_1=\mathcal{F}_2$.
We write $\gamma=\gamma_1\gamma_2$, where $\gamma_1$ is a subcurve of $\gamma$,
which lies in $Y_{l-1,1}$, has endpoints in $\mathcal{F}_1$ and $|\gamma_1|_e>0$.
If $|\gamma_2|_e=0$, then $\gamma_2$ lies in $\mathcal{F}_1$ and so $\gamma$ lies in $Y_{l-1,1}$, that contradicts
to the minimality of $l$. Hence $|\gamma_2|_e>0$.

Let $\gamma_3$ be a path in $\mathcal{F}_1$ from the terminal point of $\gamma_1$ to the initial one.
Then $|\gamma|_e>|\gamma_1|_e=|\gamma_1\gamma_3|_e$
and $|\gamma|_e>|\gamma_2|_e=|\gamma_3^{-1}\gamma_2|_e$.
Since both $\gamma_1\gamma_3$ and $\gamma_3^{-1}\gamma_2$ meet the face $\mathcal{F}_1$ and lie in $Y_l$,
they are contractible by induction.
So, $\gamma$ is contractible.

{\it Subcase 2.2.} Suppose that $\mathcal{F}_1\neq \mathcal{F}_2$.
Denote by $\widetilde{\mathcal{N}}_2$ the $\widetilde{\mathcal{N}}$-space, which is adjacent to $Y_{l-1,1}$ through the common face $\mathcal{F}_2$. Clearly
it has level $l$.
The curve $\gamma$ must leave ${\widetilde{\mathcal{N}}}_2$ through
a face $\mathcal{F}_3$. Since ${\widetilde{\mathcal{N}}}_2$ is contractible,
we may assume that $\mathcal{F}_3\neq \mathcal{F}_2$.
Continuing, we obtain that $\gamma$
passes through a cyclic sequence of subspaces  $${\widetilde{\mathcal{N}}}_1,Y_{l-1,i_1},
{\widetilde{\mathcal{N}}}_2,Y_{l-1,i_2},\dots ,{\widetilde{\mathcal {N}}}_p,Y_{l-1,i_p},$$
where we may assume that for every three consecutive subspaces $U,V,W$ the faces $U\cap V$ and $V\cap W$ are different.
So, $p\leqslant |\gamma|_e<C$, that contradicts to the construction of $Y_l$ according to Theorem~\ref{refgraph}.
$\Box$

\begin{lem}\label{loops} Every regular loop $\gamma$ in $Z$, which has length smaller than $C$, is either contractible or lies in some $Y_{1,i}$.
\end{lem}

{\it Proof.} This follows from the previous lemma in view of the facts that $Z\setminus \overset{m}{\underset{i=1}{\cup}} \widetilde{\mathcal{N}}_i$ is the disjoint union of the interiors of  $Y_{1,1},\dots, Y_{1,n}$
and that each $\widetilde{\mathcal{N}}_i$ is contractible.
$\Box$

\begin{lem}\label{nonconj} $H_2$ is not conjugate into $H_3$ in the group $G$.
\end{lem}

{\it Proof.} Assume the contrary, say $H_2^g\leqslant H_3$ for some $g\in G$,
and represent the elements $h_i^g$ and $(h_ih_j)^g$ by loops $l_i$ and $l_{ij}$ in $Z$ based at $y$.
By definition of the constant $C$, every such loop
can be freely homotoped in $Z$ to a regular loop of $e$-length smaller than $C$.
By Lemma~\ref{loops}, it can be further freely homotoped into some $Y_{1,t}$.
This means that every element $h\in \{h_i^g$, $(h_ih_j)^g: i,j=1,\dots, r\}$ can
be conjugated into some $\pi_1(Y_{1,t})=H_1^{g_t}$ by an element of $H_3$.
Then, by (3) and by Lemma~\ref{Serre3}, $H_2^g$ can be conjugated into some $H_1^{g_s}$.
This contradicts to the assumption, that $H_2$ cannot be conjugated into $H_1$ in $G$.
$\Box$

The proof of Theorem~\ref{main2} is completed.

\section{Acknowledgements} The first named author thanks the MPIM at Bonn for its support and excellent working conditions during the fall 2010, while this research was finished.






\end{document}